\documentclass[a4, 10 pt, conference]{ieeeconf}
\IEEEoverridecommandlockouts                                             % Needed to meet printer requirements.
  \usepackage{braket,amsfonts}
\usepackage{bm}
\usepackage{array}
\usepackage{placeins,graphicx,diagbox }
\usepackage{etex}
\usepackage{amsmath,amssymb}  % assumes amsmath package installed
\usepackage{pgfplots}

\usepackage{amsmath,algorithm} % assumes amsmath package installed
\usepackage{amssymb}  % assumes amsmath package installed
\usepackage{color}
\usepackage{cite}
\newcounter{algo}

\newtheorem{theorem}{Theorem}
\newtheorem{lemma}{Lemma}

\newtheorem{assumption}{Assumption}

\newtheorem{remark}{Remark}

\usepackage{algorithmic}

\usepackage{graphicx,epstopdf}

\usepackage{amsopn}

\usepackage{xspace}
\usepackage{bold-extra}
\usepackage[most]{tcolorbox}
\usepackage{placeins,diagbox }

\usepackage{wrapfig} % for wrapping the text around the fig
\usepackage{tikz,pgfplots,mathtools,float}
\usepackage{booktabs,caption}
\usepackage[flushleft]{threeparttable}
\usepackage{subfigure}
\usepackage{multirow}
\usepackage{amssymb,latexsym}  % assumes amsmath package installed
\usepackage{extarrows}
\usepackage{dsfont} % for double stroke font
\usepackage{amssymb}  % assumes amsmath package installed
\usepackage{color}
\usepackage{cite}

\colorlet{texcscolor}{blue!50!black}
\colorlet{texemcolor}{red!70!black}
\colorlet{texpreamble}{red!70!black}
\colorlet{codebackground}{black!25!white!25}

\def\Real{\mathbb{R}}

\title{Linearly Convergent  Algorithm with Variance Reduction for Distributed Stochastic Optimization }

\author{Jinlong Lei \textsuperscript{a,b},
Peng Yi \textsuperscript{a,b},
Jie Chen\textsuperscript{a,b}\quad {\it IEEE Fellow}, and
Yiguang Hong\textsuperscript{c}\quad {\it IEEE Fellow}
\thanks{\noindent\textsuperscript{a} The Department of Control Science and Engineering,
Tongji University,  Shanghai,   201804, China; }
\thanks{\noindent\textsuperscript{b}  Shanghai Institute of Intelligent Science and Technology, Tongji University,  Shanghai, 201804,  China; }
\thanks{\noindent\textsuperscript{c}    Key Laboratory of Systems and Control,
Academy of Mathematics and Systems Science, Chinese Academy of Sciences, Beijing 100190,  China;}
\thanks{ Email address: ustcleijinlong@gmail.com (J. Lei),  yipeng@amss.ac.cn (P. Yi),
 chenjie@bit.edu.cn(J. Chen) yghong@iss.ac.cn (Y.Hong).}
\thanks{ The  work of Jinlong Lei and Peng Yi was partially
supported by Key research and development projects of the Ministry of
Science and Technology of China No. 2018YFB1305304 and the National Natural Science Foundation of China under Grant No. 51475334. }
}

\begin{document}
\allowdisplaybreaks
\maketitle

%% ------------------------------------------------------------------
%% ABSTRACT
%% ------------------------------------------------------------------
\begin{abstract}
 This  paper considers a distributed   stochastic  strongly convex optimization, where  agents  over a network aim to cooperatively minimize the average of  all agents'  local cost functions.
%This  paper considers a distributed   stochastic optimization problem, where there {is a group of   agents connected over a network} and  each of them has an expectation-valued  smooth {strongly} convex cost  function. The  agents  want to collaboratively  find an optimal solution that minimizes   the average of all agents'  cost functions.
Due to the stochasticity of gradient estimation and distributedness of local objective,
  fast linearly convergent distributed algorithms have not been achieved yet.
This work proposes a  novel distributed stochastic gradient tracking algorithm with variance reduction,
where  the  local gradients are estimated by an increasing batch-size of  {sampled gradients}.
With an undirected connected communication graph   and  a geometrically increasing  batch-size,
  the {iterates}  are  shown to converge in mean to the optimal solution at a geometric rate (achieving {linear} convergence).
The iteration, communication, and oracle complexity for obtaining an  $\epsilon$-optimal solution are established as well.
Particulary, the communication complexity is  $\mathcal{O}(\ln (1/\epsilon))$ while
the oracle complexity (number of {sampled gradients}) is $\mathcal{O}(1/\epsilon^2)$, {which is of the  same order as that of centralized approaches}.
 Hence, the proposed scheme is communication-efficient  { without requiring extra  sampled gradients}.
Numerical  simulations  are given to demonstrate the theoretic results.
\end{abstract}

%\begin{IEEEkeywords}
%  proximal stochastic gradient,  randomized block-coordinate descent, nonsmooth optimization,  nonconvex stochastic optimization, variable sample-size schemes
%\end{IEEEkeywords}

\section{Introduction}
Distributed optimization %,  where agents connected over networks aim to cooperatively minimize the average of  all agents'  local cost functions,
 has  been extensively studied in recent years due to its  wide applications in
sensor networks \cite{towfic2015stability,rabbat2004distributed},
power systems \cite{yi2016initialization,yu2018economic},   distributed estimation and control   \cite{fang2018optimization,wang2018distributed,abdelatti2018cooperative}.
Various distributed optimization methods have been  developed, including
 primal domain   methods  \cite{nedic2009distributed,pu2018push},
 dual domain methods \cite{palomar2006tutorial,mota2013d}, and primal-dual domain methods  \cite{chang2014distributed,yi2015distributed,lei2016primal}.  Please refer to the survey   \cite{yang2019survey} for  more references.
%including (i)  primal domain   methods
% that  combine (sub)gradient steps with  local averaging,
% such as  distributed subgradient methods \cite{nedic2009distributed},
% distributed Nesterov gradient \cite{jakovetic2014fast},
%  and distributed gradient tracking  \cite{pu2018push}; (ii) dual domain methods employing the Lagrangian dual,
%  e.g., distributed dual decomposition \cite{palomar2006tutorial} and distributed ADMM \cite{mota2013d};
%  and (iii) primal-dual domain methods  \cite{chang2014distributed,yi2015distributed,lei2016primal}.

This paper aims to provide a fast and communication-efficient algorithm
for  distributed stochastic optimization.  %,which has been a standard model in distributed machine learning \cite{forero2010consensus}.
 We {propose a  stochastic  variant of}  the distributed gradient tracking  scheme \cite{pu2018push},
where each  agent is equipped with an auxiliary variable to  track the  dynamical average gradient in addition to the solution estimate. To achieve fast convergence and save communication cost (which is usually hundreds times of local computation cost),
   {an adaptive sampling method is incorporated into the scheme.}
The main   contributions of the paper are given as follows.
\begin{itemize}
\item We first combine the distributed stochastic gradient tracking algorithm with adaptive variance reduction,
where each agent estimates its local gradients with an increasing batchsize  of sampled  gradients.
Then each agent takes   a weighted average of   its neighbors'  estimates  and moves towards   the
negative direction of its local noisy gradient estimation.

\item When each agent's local objective function is strongly convex with  Lipschitz-continuous gradient and
the sample size for gradient estimation adaptively increases at a geometric rate,
  the  proposed   scheme  with a constant stepsize can generate geometrically/linearly convergent iterates.

\item  Furthermore, it is  shown  that  the iteration, communication, and oracle complexity
for  {each agent $i$} to obtain  an $\epsilon$-optimal solution are
$ \mathcal{O}(\ln(1/\epsilon))$,   $ \mathcal{O}(|\mathcal{N}_i|\ln(1/\epsilon))$, $\mathcal{O}(1/\epsilon^2)$,
 respectively,   where $|\mathcal{N}_i|$ denotes  the number of agent $i$'s neighbors.
Compared with existing distributed methods,
 the scheme saves the communication cost without increasing the overall sampling burden.
\end{itemize}

{\it Literature review on   distributed  stochastic optimization}.
Considerable  works have been done in distributed stochastic optimization,
e.g.,   distributed stochastic subgradient
projection algorithm  \cite{ram2010distributed},
distributed asynchronous algorithm \cite{srivastava2011distributed},  and
distributed primal-dual    method \cite{lei2018asymptotic}.
In the following, we review   some literature  on the convergence  rate and complexity analysis of distributed stochastic {strongly convex} optimization.
The work  \cite{jakovetic2018convergence}  proposed a distributed stochastic gradient method over a random network
and established  the convergence rate of    $\mathcal{O}(1/k) $   in a mean-squared sense.
A distributed stochastic mirror descent method with rate  $\mathcal{O}(\ln(k)/k) $
was given in \cite{yuan2018optimal} for  non-smooth functions,
while  a  stochastic subgradient descent with rate  $\mathcal{O}(n\sqrt{n}/k) $
 was proved in \cite{sayin2017stochastic}.
The work \cite{nedic2016stochastic} proposed a subgradient-push method  over time-varying directed graphs
and  obtained a    rate   $\mathcal{O}(\ln(k)/k) $.
A distributed stochastic  gradient tracking method
with a constant stepsize was  designed in \cite{pu2018distributed},
which  only showed that the iterates  are attracted to a {\it neighborhood}
 of the optimal solution in expectation with an exponential rate,
however, the exact convergence can not be achieved yet.

This work considers minimizing {smooth objectives} over an undirected connected network.
Instead of  decaying stepsizes in  \cite{jakovetic2018convergence,yuan2018optimal,sayin2017stochastic,nedic2016stochastic},
we adopt a constant stepsize and achieve exact and fast convergence.
%, e.g., geometric convergence rate,  endowing
%the agents with fast adaptivity in time-vary environment.
By progressively  reducing {the   variance of  gradient noises through increasing sample size},
the derived iteration complexity   matches that of centralized approaches for  deterministic optimization,
achieving superior complexity bounds than the prior works
\cite{jakovetic2018convergence,yuan2018optimal,sayin2017stochastic,nedic2016stochastic,pu2018distributed}.
Moreover,  the approach can  significantly reduce the communication rounds, {meanwhile}
the oracle complexity  can be comparable with  existing distributed stochastic gradient algorithms.
To the best of our knowledge, this is the first work that achieves a  {linear} convergence rate
for  distributed strongly convex stochastic optimization.

The paper is organized as follows. A  distributed   stochastic gradient tracking algorithm with variance reduction is proposed  in  Section II. The geometric convergence rate along with  the complexity bounds are established in Section III.
The numerical studies are presented  in Section IV, while concluding remarks are given in Section V.

 {\em Notations}.     Depending on the argument, $|\cdot|$ stands for the absolute value of a real number or the cardinality of a set.  The Euclidean norm of a vector or a matrix  is denoted as $\|\cdot\|_2$ or $\|\cdot\|$.
Let $\otimes$ denote  the Kronecker  product.
  Let  $ \mathbf{1}  $ denote   the column vectors with all entries equal to 1 and $ I_d $  denote the $d\times d$ identity matrix.
  An undirected graph  is denoted by  $\mathcal{G}=\{ \mathcal{V},\mathcal{E}\},$ where $\mathcal{V}=\{1,\dots,n\}$  is a finite set of  nodes and each   edge  $(i,j) \in \mathcal{E}$ is an unordered  pair of two distinct  nodes  $i,j$.
  A    path in $\mathcal{G}$  from $v_1$ to $v_{p }$ is a  sequence of distinct nodes, $v_1 \dots v_{p}$,
  such that  $(v_m, v_{m+1}) \in \mathcal{E}$ for all $m=1,\dots,p-1$. The graph $\mathcal{G}$ is termed {\it   connected} if for any two distinct nodes $i,j\in\mathcal{V}$,
   there is a  path between them.
   The set of node $i$'s  neighboring nodes, denoted  by $\mathcal{N}_i$,
  is defined as  $\mathcal{N}_i\triangleq \{ j\in \mathcal{V}: (i,j)\in \mathcal{E}\}$.
Define the  adjacency matrix of  graph $\mathcal{G}$ as $A=[a_{ij}]_{i,j=1}^n$, where $a_{ij}>0$  if $j\in \mathcal{N}_i$ and $a_{ij}=0$  otherwise.

\section{Problem  Statement and Algorithm Design}
In this section, we first formulate a distributed stochastic optimization problem with some assumptions.
Then we propose a  fully distributed  variable sample-size stochastic gradient tracking algorithm.

\vskip 5mm

\subsection{Problem Formulation}
We consider a network of   $n$ agents  indexed as $\mathcal{V}=\big\{1,\dots,n\big\}$, where the agents interaction   is described by an undirected     graph $\mathcal{G}=\{ \mathcal{V},\mathcal{E}\}$.
Agent $i \in \mathcal{V}$ has an expectation-valued  cost function  $f_i(x)\triangleq \mathbb{E}_{\xi_i}[h_i(x,\xi_i)]$, where  $x\in \mathbb{R}^d $,  the random vector  $\xi_i: { \Omega_i} \to \Real^{m_i}$ is
 defined on the probability space $({ \Omega_i}, {\cal F}_i, \mathbb{P})$,
and $h : \mathbb{R}^d\times \mathbb{R}^{m_i} \to \mathbb{R}$ is a scalar-valued  function. The    agents in the network need to cooperatively find an optimal solution that  minimizes the average of all agents' local cost functions, i.e.,
  \begin{equation}\label{problem1}
\begin{split}
\min_{x\in \mathbb{R}^d }  F(x )\triangleq {1\over n } \sum_{i=1}^n f_i(x).
\end{split}
\end{equation}

 We aim to design  a   distributed algorithm to drive all agents'  iterates to the optimal solution, explore  its convergence  rate,
 and establish  the  complexity bounds  for obtaining an optimal solution with a prescribed accuracy.
Below are the assumptions on the communication graph and cost functions.
 \begin{assumption}\label{ass-graph} The  undirected  graph  $\mathcal{G}$  is  connected, and its adjacency matrix $A=[a_{ij}]_{i,j=1}^n$ is   symmetric   with  the weights  $a_{ij}$   satisfying the following condition:
 \begin{align}
  \sum_{i=1}^n a_{ij}=1, \quad \forall j\in\mathcal{V} \label{adj-a}   .
 \end{align}
  \end{assumption}

With Assumption \ref{ass-graph},  the adjacency matrix  $A$ of the connected communication  graph is doubly stochastic.
 It has been shown  in \cite{xiao2004fast}  that the   spectral radius     $\sigma_A$   of
$A-\mathbf{1}\mathbf{1}^T/n$  satisfies  $\sigma_A\in (0,1)$.
%and specially      $\sigma_A=\max\{\lambda_2, \lambda_n\}$  .
% Then we have  the following relation, which will play a important  rule in the algorithm analysis part.
%  \begin{align}\label{bd-consensus}
%&  \| (A  \otimes I_d ) x-   (\mathbf{1} \otimes I_d)  \bar{x}\|  = \| (A-\mathbf{1}\mathbf{1}^T/n)\otimes I_d  \big(x- (\mathbf{1} \otimes I_d)  \bar{x} \big)  \|  \notag\\
%&   \leq   \| A-\mathbf{1}\mathbf{1}^T/n \| \|   x- (\mathbf{1} \otimes I_d)  \bar{x} \|   \leq     \sigma_A \|   x - (\mathbf{1} \otimes I_d)  \bar{x}  \| .
%\end{align}

 \begin{assumption}\label{ass-fun}
    For each agent $i\in \mathcal{V},$ its cost function $f_i(x)$ is $\eta$-strongly convex
    and its gradient function is  $L$-Lipschitz continuous, i.e.,  for any $x_1,x_2 \in \mathbb{R}^d$:
\\ (i) $(\nabla f_i(x_1)-  \nabla f_i(x_2) )^T(x_1-x_2) \geq \eta \| x_1-x_2\|^2,$   \\ (ii)  $  \| \nabla f_i(x_1)-  \nabla f_i(x_2) \|  \leq L \| x_1-x_2\|.$   \end{assumption}

By Assumption \ref{ass-fun} and  definition $F(x )= {1\over n } \sum_{i=1}^n f_i(x)$,
 $F(x)$ is $\eta$-strongly convex and its gradient function  is $L$-Lipschitz continuous.  Then problem \eqref{problem1} has a unique optimal solution, denoted by  $x^*$.
 Hence, $\nabla F(x^*)=0$  by the first-order  optimality condition.

Suppose  there exists a {\em stochastic first-order oracle} for each agent $i \in \mathcal{V}$
 such that for   any  given $x,\xi$, a  sampled gradient $\nabla h_i(x,\xi)$ is returned, which   is  an unbiased estimator of  $\nabla f_i(x)$ with bounded second-order moment. Here is the assumption on the stochastic first-order oracle.

  \begin{assumption}\label{ass-noise}  There exists a constant $\nu>0$ such that    the following holds
  for each $i\in \mathcal{V}$  and any given $x\in \mathbb{R}^d$,
  \begin{align*}
&  \mathbb{E}_{\xi_i}[\nabla h_i(x,\xi_i)]=\nabla f_i(x),  {~\rm and ~} \\& \mathbb{E}_{\xi_i}[\| \nabla h_i(x,\xi_i)-\nabla f_i(x)\|^2]  \leq \nu^2 .
  \end{align*}
    \end{assumption}

\vskip 5mm

\subsection{A Distributed   Stochastic Gradient Tracking  Algorithm with Variance Reduction }

The discrete time is slotted at $k=0,1,2,\dots$. Each agent $i$ at time $k$ maintains two estimates $x_i(k)$ and   $y_i(k)$,
where  $x_i(k)$ and $y_i(k)$ are    used  to estimate the optimal solution and to  track  the average gradient, respectively.
 Since  the  exact gradient of each expectation-valued cost function $f_i(x)$ is unavailable, we  approximate it by  averaging through  a variable batch-size  of  sampled gradients:  \begin{align}\label{est-grad}
 \tilde{g}_i(x_i(k))={1\over N(k)} \sum_{p=1}^{N(k)}\nabla  h_i(x_i(k),\xi^p_i(k)), \quad \forall k\geq 0,
 \end{align}
where   $N(k)$   is the number of sampled gradients  utilized at  time   $k $ and    the  samples $\{ \xi^p_i(k)\}_{p=1}^{N(k)}$ are   randomly and independently  generated  from the  probability space $({ \Omega_i}, {\cal F}_i, \mathbb{P})$.
The gradient estimate  given by \eqref{est-grad} is an unbiased estimate of the exact gradient,
and the  variance of the gradient noise will be progressively reduced by increasing  the batch-size.
 By  combining the distributed gradient tracking scheme   \cite{pu2018push}  with a variance reduction scheme, we  obtain  Algorithm~\ref{alg-sto-gradient_track}.
We will specify the selection of the constant steplength  $\alpha$ and the batch-size $N(k)$ upon convergence analysis.

\begin{algorithm}
\caption{A distributed variable sample-size stochastic gradient tracking algorithm}\label{alg-sto-gradient_track}
{\em Initialization}: Set $k:=0$. For any $i = 1, \hdots, n$, let  $y_i(0)= \tilde{g}_i(x_i(0))$ with  arbitrary initial  $ x_i(0) \in  \mathbb{R}^{d} $.

{\em Iterate until convergence.}

Each agent  $i=1,\cdots, n$ updates  its estimates as follows:
\begin{subequations}
\begin{align}
 x_i(k+1)& =\sum_{j\in \mathcal{N}_i}  a_{ij} x_j(k)-\alpha y_i(k), \label{push}\\
 y_i(k+1)&= \sum_{j\in \mathcal{N}_i} a_{ij}  y_j(k) +\tilde{g}_i(x_i(k+1))-  \tilde{g}_i(x_i(k)), \label{pull}
  \end{align}
\end{subequations}
 where   $\alpha>0$ is the steplength and  $ \tilde{g}_i(x_i(k))$ is  given  in \eqref{est-grad}.
   \end{algorithm}

Note that for each agent $i\in \mathcal{V},$ the  implementation of   Eqn. \eqref{push}  requires  its neighbors'  estimates of the optimal solution  $\{ x_j(k)\}_{ j\in \mathcal{N}_i}$,
 while the update of $y_i(k+1) $    characterized by Eqn.  \eqref{pull}  uses its  local gradient estimate  as well as  its neighbors'  information $\{ y_j(k)\}_{ j\in \mathcal{N}_i}$ to  asymptotically track   the dynamical average  gradient across the network.
  Therefore,  Algorithm~\ref{alg-sto-gradient_track} is a fully distributed algorithm since the update of each agent merely uses its local data and  its neighboring   information.

\section{Convergence Analysis}

In this section, we show the geometric convergence rate  for Algorithm~\ref{alg-sto-gradient_track} when the batchsize is increased at  a geometric rate,
and  establish   the complexity bounds for obtaining an $\epsilon$-optimal solution.

\vskip 5mm
\subsection{Preliminary Lemma}

Define  the  gradient observation noise  as follows:
\begin{equation}\label{def-w}
  w_i(k)\triangleq  \tilde{g}_i(x_i(k)) -\nabla f_i(x_i(k)).
\end{equation}    We further define
\begin{equation}\label{def-fw}
\begin{array}{ll}
   & x (k) \triangleq \big( x_1(k) ^T, \cdots, x_n(k) ^T \big)^T,   \\& y(k ) \triangleq \big(y_1(k )^T, \cdots, y_n(k)^T \big)^T  ,    \\&  \nabla (k) \triangleq \big(\nabla f_1(x_1(k)) ^T, \cdots, \nabla f_n(x_n(k)) ^T \big)^T,
\\&   w(k ) \triangleq \big(w_1(k )^T, \cdots, w_n(k)^T \big)^T  .
\end{array}
\end{equation}
Then  Algorithm \ref{alg-sto-gradient_track} is written in a compact form:
 \begin{subequations}
\begin{align}
 x (k+1)& =(A  \otimes I_d ) x(k)-\alpha y(k)  \label{push2}\\
 y(k+1)&= (A  \otimes I_d )   y(k) +\nabla (k+1)   \notag\\&\quad +w(k+1)-\nabla (k)-w(k ).\label{pull2}
  \end{align}
\end{subequations}

Define the average of  agents'  estimates of the optimal  solution and  the averaged gradient across the network  as follows for any  $k\geq 0$:
\begin{align}\label{def-xy}
\bar{x}(k)={1\over n} \sum_{i=1}^n x_i(k)  {\rm~and}~\bar{y}(k)={1\over n} \sum_{i=1}^n y_i(k) .
\end{align}
We  start  to analyze  the algorithm performance  by   characterizing the   interactions  among the three   error sequences: (i)   distance from the average estimate to   the optimal solution $\| \bar{x}(k) -x^*\|;$
(ii) consensus error $  \| x(k) - (\mathbf{1} \otimes I_d)\bar{x}(k)  \| $;
and  (iii) consensus  error  of the gradient trackers   $\| y(k) - (\mathbf{1} \otimes I_d)  \bar{y}(k)\|$.
We     bound the three  error sequences in terms of the linear combinations of their past values in the following lemma,
of which  the proofs can be found  in Appendix.

\begin{lemma}\label{lem1}  Suppose Assumptions \ref{ass-graph} and \ref{ass-fun} hold.   Let the sequences  $\{x(k)\}$ and $\{ y(k)\}$ be generated by  Algorithm \ref{alg-sto-gradient_track} with $0<\alpha \leq {2\over \eta+L}.$
 Define  the following vector   and matrix:
\begin{equation} \label{def-z}
\begin{split} z(k) & \triangleq
\begin{pmatrix}
\| \bar{x}(k) -  x^*  \| \\
\| x(k) - (\mathbf{1} \otimes I_d)\bar{x}(k)\|  \\
 \| y(k) - (\mathbf{1} \otimes I_d)  \bar{y}(k)\|
\end{pmatrix} ,  \\
J(\alpha) & \triangleq
\begin{pmatrix}
1-\alpha L  &  \alpha {L \over \sqrt{n}}  & {0}  \\
 {0} &   \sigma_A & \alpha  \\
 \alpha \sqrt{ n} L^2 &L    \|  A-I_n \|  + \alpha L^2  &\sigma_A+\alpha L
 \end{pmatrix} .
 \end{split}
 \end{equation}
 Then   the following inequalities hold for  any $k\geq 0$:
\begin{equation}\label{recursion-z0}
\begin{split}
& z(k+1) \leq  J(\alpha) z(k)    \\&+\begin{pmatrix}
 {\alpha \over n} \sum_{i=1}^n\left\| w_i(k)\right \|\\
{0} \\ \| w(k+1)-w(k ) \|  +{\alpha   L \over \sqrt{n}} \sum_{i=1}^n \left \|   w_i(k) \right \|
\end{pmatrix} .
\end{split}
\end{equation}
\end{lemma}

\vskip 5mm

\subsection{Linear convergence rate analysis}

We now give the linear convergence rate result.
\vskip 5mm
\begin{theorem}  \label{thm1}  Suppose Assumptions \ref{ass-graph},  \ref{ass-fun},  and \ref{ass-noise} hold.
 Let  $\{x(k)\}$ and $\{ y(k)\}$ be generated by  Algorithm \ref{alg-sto-gradient_track} with $0<\alpha \leq {2\over \eta+L}.$
 Set $N(k)=\lceil q^{-2k} \rceil$ for some $q\in (0,1). $
  Take $\alpha$ such that the spectral radius of the matrix $   J(\alpha)$, denoted by $\rho(J(\alpha)) $,  is strictly smaller than 1.
  Then  $z(k)$  converges to zero in mean at a geometric rate, that is,
 \begin{subequations} \label{linear-rate}
\begin{align}
\mathbb{E}[\| z(k)\|] &  \leq   \rho(J(\alpha))^k  \| \mathbb{E}[z(0) ]\|  +   {C \over q-\rho(J(\alpha))} q^k ,   \notag
\\&   \qquad \qquad  { \rm when ~}  q>\rho(J(\alpha)) ,    \label{linear-rate_a}
\\  \mathbb{E}[\| z(k)\|]  & \leq  \rho(J(\alpha))^k  \| \mathbb{E}[z(0) ]\|  +   {C \over  \rho(J(\alpha)) -q} \rho(J(\alpha))^k   , \notag
\\& \qquad  \qquad{~\rm when ~}  q< \rho(J(\alpha)),  \label{linear-rate_b}
  \end{align}
\end{subequations}
where  $C\triangleq  \nu \sqrt{ \alpha^2 +  n  (   1+q  + \alpha   L   )^2}.$
\end{theorem}
{\bf Proof.}  We first  split the matrix $ J(\alpha) $   into the sum of a fixed matrix and
another perturbation matrix as a function $\alpha$:
\begin{equation}
\begin{array}{ll}
   J(\alpha) & \triangleq  \begin{pmatrix}
1&  0 & 0\\ 0&   \sigma_A & 0 \\
0 &L    \|  A-I_n \|    &\sigma_A  \end{pmatrix}   +\alpha \begin{pmatrix}
-L &    {L \over \sqrt{n}}  & 0\\ 0& 0 & 1  \\
  \sqrt{ n} L^2 &  L^2  &  L \end{pmatrix}  \\
  &
  := J(0)+\alpha E.
  \end{array}\nonumber
  \end{equation}

Because  $\sigma_A\in(0,1),$ the spectral radius of $J(0)$ is 1   and  the  corresponding right and
left eigenvector  to the eigenvalue of 1 of  $J(0)$  is $(1,0,0)^T$.  Because the eigenvalues of a matrix  are  a  continuous function of its entries,
we are able to choose some sufficiently small $\alpha$ such that the spectral radius of $   J(\alpha)$ is strictly smaller than 1  (see \cite[Theorem 1]{xin2018linear} for a more detailed discussion).

Define $\mathcal{F}_k\triangleq \sigma \left \{ x(0),   \{ \xi^p_i(t)\}_{p=1}^{N(t)}, i\in \mathcal{V},  0\leq t \leq k-1  \right \}$.
Thus,  $x(k)$ produced    by  Algorithm \ref{alg-sto-gradient_track}  is adapted to $\mathcal{F}_k $.  Then by   \eqref{est-grad},  \eqref{def-fw}, and  the fact that  the   samples $\{ \xi^p_i(k)\}_{p=1}^{N(k)}$ are    independent, we obtain that
\begin{equation*}
\begin{split}
& \mathbb{E}[\| w_i(k) \|^2 |F_k ]
%\\&={1 \over N(k)^2}\mathbb{E} \Big[  \| \sum_{p=1} ^{N(k) }   (\nabla  h_i(x_i(k),\xi^p_i(k))- \nabla  f_i(x_i(k))  ) \|^2 \big|F_k  \Big]
\\& ={1\over N(k)^2} \sum_{p=1} ^{N(k) } \mathbb{E}   [  \|   \nabla  h_i(x_i(k),\xi^p_i(k))- \nabla  f_i(x_i(k))    \|^2 |F_k   ] .
\end{split}
\end{equation*}
Then by Assumption \ref{ass-noise},  the following holds for each  $i\in \mathcal{V}$:
\begin{equation*}
\begin{split}
& \mathbb{E}[\| w_i(k) \|^2] \leq {\nu^2 \over N(k)}, \quad \forall k\geq 0.
\end{split}
\end{equation*}
Therefore,  from $N(k)=\lceil q^{-2k} \rceil$  and  the relation  $\mathbb{E}[\|x\|] \leq\sqrt{\mathbb{E}[\|x\|^2] } $,  we have for any $k\geq 0$,
\begin{equation*}
\begin{split}
\mathbb{E}[\| w_i(k) \|]&  \leq \sqrt{\mathbb{E}[\|w_i(k)\|^2] } \leq {\nu  \over  \sqrt{N(k)}} =\nu q^k,
  \\ \mathbb{E}[\| w (k) \|] &\leq\sqrt{\mathbb{E}[\|w(k)\|^2] }=\sqrt{\sum_{i=1}^n \mathbb{E}[\|w_i(k)\|^2] } \leq \sqrt{n} \nu q^k.
\end{split}
\end{equation*}
Then by taking expectations on both sides of Eqn. \eqref{recursion-z0} and using the triangle equality, we obtain the following  entry-wise  linear matrix  inequality:
\begin{align*}%\label{recursion-z1}
&\mathbb{E}[z(k+1) ] \leq  J(\alpha) \mathbb{E}[z(k) ] +\begin{pmatrix}
  \alpha \nu  \\ 0 \\  \sqrt{n} \nu (1+q) + \alpha   L   \sqrt{n} \nu
\end{pmatrix}  q^k \notag  \\
&\leq  J(\alpha) ^{k+1}\mathbb{E}[z(0) ] + \sum_{p=0}^k  J(\alpha)^p  \begin{pmatrix}
  \alpha \nu  \\ 0 \\  \sqrt{n} \nu (1+q + \alpha   L )
\end{pmatrix}  q^{k-p} .
\end{align*}
Therefore,  we can obtain the following bound for any $k\geq 1$:
\begin{equation}\label{bd-errz}
\begin{split}
& \| \mathbb{E}[z(k ) ] \|
\leq   \| J(\alpha) \| ^{k } \| \mathbb{E}[z(0) ]\|
\\&+ \sum_{p=0}^{k-1}   \| J(\alpha) \| ^p  \left \|   \begin{pmatrix}
  \alpha \nu  \\ 0 \\  \sqrt{n} \nu (1+q + \alpha   L )
\end{pmatrix}  \right \|   q^{k-1-p}  \\&  \leq
\rho(J(\alpha))^k  \| \mathbb{E}[z(0) ]\|  +  C  \sum_{p=0}^{k-1}  \rho(J(\alpha)) ^p    q^{k-1-p}  .
\end{split}
\end{equation}
  Note that  for any $\rho<q$:
\begin{align*}
\sum_{p=0}^{k-1}  \rho ^p    q^{k-1-p}  =q^{k-1}  \sum_{p=0}^{k-1}  (\rho/q) ^p  \leq {q^{k-1} \over 1-\rho/q} ={q^k \over q-\rho},
\end{align*}
while for  any $\rho>q$: $
\sum_{p=0}^{k-1}  \rho ^p    q^{k-1-p}     \leq   {\rho^k \over  \rho-q}.$
Combining with  \eqref{bd-errz}, we prove the geometric rate \eqref{linear-rate}.
\hfill $\Box$

 It is noticed from Algorithm \ref{alg-sto-gradient_track}    that  all agents   use an identical steplength $\alpha$,
 which  may require additional coordination among the agents before running     the   algorithm.
Recently, techniques utilizing uncoordinated steplengths have been proposed in \cite{nedic2017achieving}. How to incorporate such a scheme with the
variance reduced method remains our future work.
Besides,  in Theorem \ref{thm1}, $\alpha$  is chosen  to be sufficiently small such   that  $\rho(J(\alpha)) <1. $ This is merely
a sufficient condition for guaranteeing  linear convergence,   and the necessary condition on the steplegnth $\alpha$ remains an open problem.

\begin{remark}
Theorem \ref{thm1}  implies that if  the number of sampled gradients is increased at a geometric rate $\lceil q^{-2k} \rceil$ with $q\in (0,1)$,
the   expectation valued  error sequences   $ \mathbb{E}[\| \bar{x}(k) -x^*\| ],$ $ \mathbb{E}[ \| x(k) - (\mathbf{1} \otimes I_d)\bar{x}(k)  \| ]$,
and   $ \mathbb{E}[\| y(k) - (\mathbf{1} \otimes I_d)  \bar{y}(k)\|]$    converge  to zero  at a geometric rate of $\max\{ \rho(J(\alpha))^k,q^k\}$.
When $0< q<\rho(J(\alpha))$,  the  geometric     rate  $\mathcal{O}( \rho(J(\alpha))^k)$ in the  deterministic regimes might  be  recovered.
\end{remark}

\subsection{Complexity Analysis}

Based on the geometric convergence rate established in Theorem \ref{thm1},
we are able  to  establish the complexity bounds
for obtaining an  $\epsilon$-optimal solution satisfying   $ \mathbb{E}[\| z \|]  \leq \epsilon.$
 The   iteration  complexity is defined as  $K(\epsilon)$  such that $ \mathbb{E}[\| z(k )\|]  \leq \epsilon $ for any $ k \geq K(\epsilon)$.
With the updates in  Algorithm \ref{alg-sto-gradient_track},
  agent  $i$  requires $2 | \mathcal{N}_i| $ rounds of communications to obtain its neighbors' information $x_j(k), y_j(k), j\in\mathcal{N}_i.$
Thus,   the  communication complexity of agent $i$  to obtain  an   $\epsilon$-optimal solution is
 $2 | \mathcal{N}_i| K(\epsilon)$.
Agent $i$'s  oracle complexity, denoted by $O(\epsilon) $, is measured by the number of sampled gradients
 for deriving an   $\epsilon$-optimal solution, and can be computed as   $ \sum_{k=0}^{ K(\epsilon)} N(k).$
The following theorem gives the complexity bounds.

\begin{theorem}  \label{thm2}   Suppose Assumptions \ref{ass-graph},  \ref{ass-fun}, and \ref{ass-noise}  hold.
Consider  Algorithm \ref{alg-sto-gradient_track} with $0<\alpha \leq {2\over \eta+L} $  and  $N(k)=\lceil q^{-2k} \rceil$ for some $q\in (0,1). $
Take $\alpha $ such that $\rho(J(\alpha))<1,$  then the
iteration, communication, and oracle complexity  required by agent $i$  to obtain  an $\epsilon$-optimal solution  are  $ K(\epsilon)$,    $2 | \mathcal{N}_i| K(\epsilon)$, and $O(\epsilon) $, respectively,
 where  $ K(\epsilon)$ and   $O(\epsilon) $ are given as follows:
 \begin{subequations}
\begin{align}
  K(\epsilon) & =\begin{cases}  &
  \ln  \left( { \| \mathbb{E}[z(0) ]\|  +    C /(q-\rho(J(\alpha))) \over \epsilon}\right)  { 1\over \ln(1/q) } ,  \\&  \qquad \qquad  {~\rm when ~}  \rho(J(\alpha))<q,
\\ & \ln  \left( { \| \mathbb{E}[z(0) ]\|  +    C /(q-\rho(J(\alpha))) \over \epsilon}\right){ 1\over \ln(1/\rho(J(\alpha))) }    ,  \\&\qquad  {~\rm when ~}  \rho(J(\alpha))>q,
\end{cases} \label{iter_complexity} \\
O(\epsilon) & =\begin{cases}  &
{1 \over \epsilon^2(  1-q^2 )  }  \left(   \| \mathbb{E}[z(0) ]\|  +  { C \over  q-\rho(J(\alpha)) } \right)^2,
\\&   \qquad  \qquad  {~~\rm when ~}  \rho(J(\alpha))<q,
\\ &     {1\over 1-q^2}  \left( { \| \mathbb{E}[z(0) ]\|  +    C /( \rho(J(\alpha))-q) \over \epsilon}\right)^{ 2\ln(1/q) \over \ln(1/\rho(J(\alpha))) } ,
\\&  \qquad { \rm when ~}  \rho(J(\alpha))>q.
\end{cases}
\end{align}
\end{subequations}
\end{theorem}
{\bf Proof.}   We prove this theorem by considering the  two cases with  $\rho(J(\alpha))<q $ and $\rho(J(\alpha))>q$, respectively.

\noindent Case  (i).  $\rho(J(\alpha))<q$.
With Eqn. \eqref{linear-rate_a}, for any $  k\geq  K_1(\epsilon)\triangleq { 1\over \ln(1/q) }  \ln  \big( { \| \mathbb{E}[z(0) ]\|  +    C /(q-\rho(J(\alpha))) \over \epsilon}\big), $ we have\
\begin{equation}
\mathbb{E}[\| z(k )\|]  \leq  \left(   \| \mathbb{E}[z(0) ]\|  +   {C \over q-\rho(J(\alpha))} \right) q^k \leq \epsilon . \nonumber
\end{equation}
This allows us to bound agent $i$'s oracle complexity   by
\begin{align*}
&O_1(\epsilon) =\sum_{k=0}^{K_1(\epsilon)} N(k)=\sum_{k=0}^{K_1(\epsilon)}  q^{-2k}
\\&= {q^{-2(K_1(\epsilon)+1)}-1 \over q^{-2}-1} \leq {q^{-2} \over q^{-2}-1} {q^{-2K_1(\epsilon)}}
\\& \leq  {1\over 1-q^2} q^{-2 { 1\over \ln(1/q) }  \ln  \left( { \| \mathbb{E}[z(0) ]\|  +    C /(q-\rho(J(\alpha))) \over \epsilon}\right) }
\\&=  {1\over 1-q^2}  e^{\ln(q^{-2 }){ 1\over \ln(1/q) }  \ln  \left( { \| \mathbb{E}[z(0) ]\|  +    C /(q-\rho(J(\alpha))) \over \epsilon}\right) }
 \\&=  {1\over 1-q^2}  \left( { \| \mathbb{E}[z(0) ]\|  +    C /(q-\rho(J(\alpha))) \over \epsilon}\right)^2.
\end{align*}

\noindent  Case (ii). When   $\rho(J(\alpha))>q$,  by defining $K_2(\epsilon) \triangleq { 1\over \ln(1/\rho(J(\alpha))) }  \ln  \big( { \| \mathbb{E}[z(0) ]\|  +    C /(q-\rho(J(\alpha))) \over \epsilon}\big)  $,   from   Eqn. \eqref{linear-rate_b}  we have  for any $k\geq  K_2(\epsilon) :$
\begin{align*}
\mathbb{E}[\| z(k )\|]  & \leq  \left(   \| \mathbb{E}[z(0) ]\|  +   {C \over \rho(J(\alpha)) -q} \right)\rho(J(\alpha))^k   \leq \epsilon ,
\end{align*}
which allows us to bound agent $i$'s oracle complexity   by
\begin{align*}
&O_2(\epsilon) =\sum_{k=0}^{K_2(\epsilon)} N(k)=\sum_{k=0}^{K_2(\epsilon)}  q^{-2k}
\\& = {q^{-2(K_2(\epsilon)}+1)-1 \over q^{-2}-1} \leq  {q^{-2} \over q^{-2}-1} q^{-2 K_2(\epsilon)}
 \\& \leq   {1\over 1-q^2}  q^{-2 { 1\over \ln(1/\rho(J(\alpha))) }  \ln  \left( { \| \mathbb{E}[z(0) ]\|  +    C /( \rho(J(\alpha))-q) \over \epsilon}\right) }
\\&=    {1\over 1-q^2} e^{\ln(q^{-2 }){ 1\over \ln(1/\rho(J(\alpha))) }  \ln  \left( { \| \mathbb{E}[z(0) ]\|  +    C /(\rho(J(\alpha))-q) \over \epsilon}\right) }
 \\&=  {1\over 1-q^2}    \left( { \| \mathbb{E}[z(0) ]\|  +    C /( \rho(J(\alpha))-q) \over \epsilon}\right)^{2\ln(1/q) \over \ln(1/\rho(J(\alpha))) }.
\end{align*}
By combining Cases (i) and (ii), we  complete  the proof.
\hfill $\Box$

\begin{remark}  Theorem \ref{thm2}  shows that  when the bacthsize increases at a geometric rate, the
iteration complexity  required by agent $i$ to obtain an $\epsilon$-optimal solution  is  $ \mathcal{O}(\ln(1/\epsilon))$,   which is an {\em optimal}  bound for strongly convex optimization.
 Moreover, the  number of communication rounds required by agent $i$ is $ \mathcal{O}(|\mathcal{N}_i|\ln(1/\epsilon))$, which is  proportional to the number of  its   neighboring agents.
In terms of the oracle complexity, the {\em optimal}  bound $\mathcal{O}(1/\epsilon^2)$ is achieved when  the adaptive parameter $q$ satisfies $q\in (\rho(J(\alpha)),1),  $
 while  for the case with  $q\in (0,\rho(J(\alpha))), $  the suboptimal bound $  (1/ \epsilon)^{ 2\ln(1/q) \over \ln(1/\rho(J(\alpha))) }$ is  obtained because ${ \ln(1/q) \over \ln(1/\rho(J(\alpha))) }>1.$ Therefore, the communication cost is saved without increasing the sample burden.
\end{remark}

 \section{Numerical Simulations}

   In this section, we apply  Algorithm \ref{alg-sto-gradient_track}  to a  distributed parameter estimation problem   \cite{towfic2015stability}. Consider a network of  $n$ spatial
   sensors that aim to  estimate an  unknown  $d $-dimensional parameter $x^*$   in a distributed manner.
   Each sensor  $i$ collects a set of scalar    measurements   $ \{ d_{i,p} \}_{p \geq 1}$ generated by  the  following linear regression model corrupted by observation noises:
 \begin{equation}
    d_{i,p}=u_{i,p}^Tx^* +\nu_{i,p},\nonumber
 \end{equation}
  where   $u_{i,p} \in  \mathbb{R}^{d}$ is  the regression vector  accessible to  agent $i$
  and  $\nu_{i,p}  \in \mathbb{R} $ is a zero-mean random noise.

Suppose that $\{u_{i,p}\}$ and $\{\nu_{i,p}\}$ are mutually independent     Gaussian sequences with  distributions $N(\mathbf{0}, R_{u,i})$ and    $N(0,\sigma_{i,\nu}^2)$, respectively.
Then the distributed parameter estimation can be solved with a distributed stochastic  quadratic optimization problem:
 \begin{equation}\label{filter1}
  \min_{x\in \mathbb{R}^d}
  ~  {1\over n}\sum_{i=1}^n f_i(x) := \mathbb{E}   \big[\|  d_{i,p}-u_{i,p}^Tx\| ^2\big].
  \end{equation}
   Thus, $f_i(x)=(x-x^*)^T R_{u,i}(x-x^*)+\sigma_{i,\nu}^2 $ and $\nabla f_i(x)= R_{u,i}(x-x^*).$
Assume that the covariance $ R_{u,i}$  is positive definite, then   $x^*$ is the unique optimal solution to  \eqref{filter1}.
By using  the observed  regressor   $u_{i,p}$ and  the corresponding measurement $d_{i,p}$,  a  noisy sample of the  exact gradient $\nabla f_i(x) $ is $ u_{i,p}u_{i,p}^T x-d_{i,p} u_{i,p}$, satisfying Assumption \ref{ass-noise}.

\begin{figure}[!htb]
 \includegraphics[width=3.5 in]{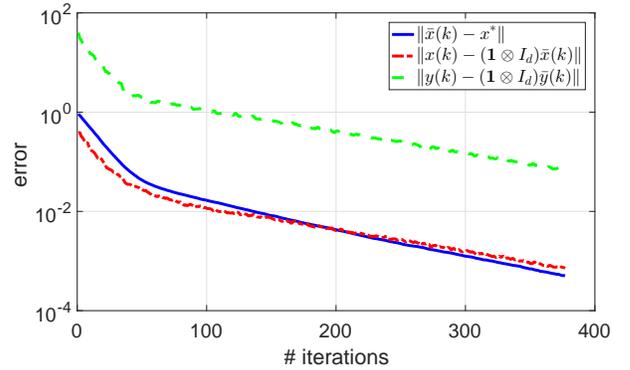}
 \caption{ Geometric Rate of Convergence }     \label{GOne}
\end{figure}

\begin{figure}[!htb]
 \includegraphics[width=3.5in]{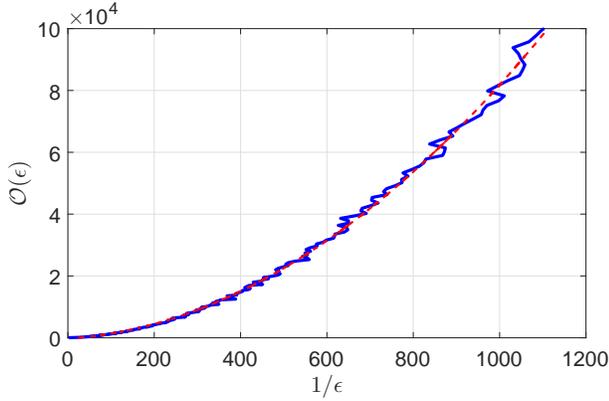}
 \caption{ Empirical oracle complexity and its quadratic fitting}     \label{GThree}
\end{figure}

In the   experiment, we set $x^*=\mathbf{1}/\sqrt{d},$   $d=5,$ and $n=10.$
We  randomly generate an  undirected and connected network,  where any two  distinct agents are  linked with probability $0.3.$
The adjacency matrix is constructed based on the Metropolis rule \cite{sayed2014adaptive}.  We now run Algorithm \ref{alg-sto-gradient_track} with $\alpha=0.01$ and $N_k=\lceil 0.98^{-k}\rceil $ and examine the   empirical  rate  of convergence  and oracle complexity,
where the empirical mean is  calculated by averaging across 50 sample paths. The convergence rate
% of   $\mathbb{E}[\| \bar{x}(k) -x^*\|] $,   $ \mathbb{E}[  \| x(k) - (\mathbf{1} \otimes I_d)\bar{x}(k)  \|] $, and  $\| y(k) - (\mathbf{1} \otimes I_d)  \bar{y}(k)\|$
  is   shown in Fig.  \ref{GOne},  which demonstrates   that the iterates
$\{x_k\}$   generated by   Algorithm  \ref{alg-sto-gradient_track} converge  in mean to the  true parameter $x^*$ at a linear rate.  Furthermore, the   relation between  $\epsilon$ and $O(\epsilon)$  is shown in   Fig.  \ref{GThree}
with the blue solid curve representing the empirical data and the  red    dashed curve denoting its quadratic fitting,
where  $O(\epsilon)$ denotes  the   number of    sampled gradients required   to make
$$\mathbb{E}\left[\left\| \begin{pmatrix}  \bar{x}(k) -x^* \\ x(k) - (\mathbf{1} \otimes I_d)\bar{x}(k)  \end{pmatrix} \right \| \right]   <\epsilon.$$ Fig.  \ref{GThree}  implies that the empirical
oracle complexity fits well with the established theoretical bound $\mathcal{O}(1/\epsilon^2)$.

\begin{figure}[!htb]
 \includegraphics[width=3.5in]{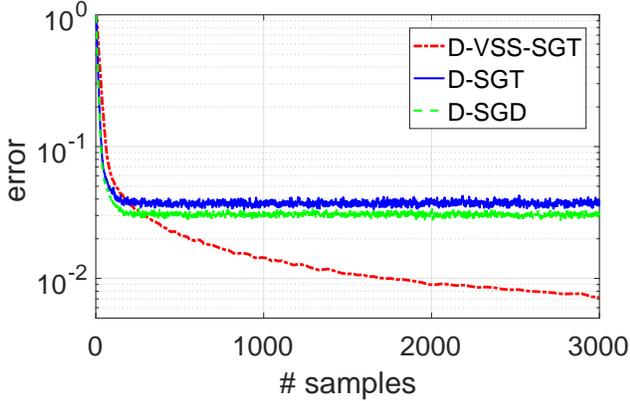}
 \caption{Comparison with  D-SGD and D-SGT}     \label{Giter}
\end{figure}

   We then  compare Algorithm  \ref{alg-sto-gradient_track}, abbreviated as D-VSS-SGT,  with the  distributed stochastic gradient  descent (D-SGD) \cite{ram2010distributed} and the distributed stochastic gradient tracking  (D-SGT) \cite{pu2018distributed}.
We set   the   consant steplength as $\alpha=0.01$ in the three  schemes,    $N_k=\lceil 0.98^{-k}\rceil $  in Algorithm \ref{alg-sto-gradient_track}, and terminate them when the  total number of  sampled gradients  utilized reaches 3000.
The empirical error  $\mathbb{E}\left[\left\| \begin{pmatrix}  \bar{x}(k) -x^* \\ x(k) - (\mathbf{1} \otimes I_d)\bar{x}(k)  \end{pmatrix} \right \| \right]  $ of the three algorithms  vs the number of sampled gradients   is given in Figure \ref{Giter}.
It can be seen that the iterates of  D-SGD  and D-SGT  ceased at a neighborhood of the   true parameter $x^*$,
while  the iterate generated by Algorithm  \ref{alg-sto-gradient_track}  converges  to  the true value  $x^*$ at a faster convergence speed.

\section{Conclusions}

We  designed  a  novel  distributed variance reduced  stochastic gradient tracking algorithm for strongly convex stochastic optimization over networks.
 We proved that with  a  suitably selected constant steplength,
the iterates converge  in mean to the optimal solution at a geometric rate when the batch-size is increased geometrically.
We further establish the   complexity bounds for obtaining an $\epsilon$-optimal solution, where the iteration complexity $ \mathcal{O}(\ln(1/\epsilon))$ matches the optimal bound of  centralized approaches in the deterministic regimes,
 the communication complexity  is significantly reduced    to  $ \mathcal{O}(|\mathcal{N}_i|\ln(1/\epsilon))$,
and the oracle complexity  $\mathcal{O}(1/\epsilon^2)$ is comparable with  the standard stochastic gradient descent algorithm.
In future, we will consider   asynchronous approaches with  agent-specific  stepsize and batchsize, and contend with the directed or switching graphs.

In future, we will consider   asynchronous approaches with  agent-specific  stepsize and batchsize, and contend with the directed or switching graphs.
It is also worthwhile investigating the communication and oracle complexity for distributed stochastic optimization with other variance reduction methods, like \cite{johnson2013accelerating},\cite{defazio2014saga} and \cite{roux2012stochastic}. The extension of the current algorithm to non-convex/non-smooth distributed stochastic optimization
is also a promising research direction.

\appendix
\section{Proof of Lemma \ref{lem1}}

{\bf Proof of Lemma \ref{lem1}.}
By  multiplying both sides of   Eqn. \eqref{push2} with ${ (\mathbf{1}  \otimes I_d ) \over n} $ from the left
and using Assumption \ref{ass-graph},  we   obtain  that
\begin{align}\label{recur-barx1}
\bar{x}(k+1)=\bar{x}(k) -\alpha\bar{y}(k), \quad \forall k\geq 0.
\end{align}
Also, by using  \eqref{pull} and Assumption \ref{ass-graph},  there holds
\begin{equation}\label{recur-bary0}
\begin{split}
\bar{y}(k+1) & =\bar{y}(k) +{1\over n} \sum_{i=1}^n \tilde{g}_i(x_i(k+1))    \\& \quad
-  {1\over n} \sum_{i=1}^n \tilde{g}_i(x_i(k)) , \quad \forall k\geq 0.
\end{split}
\end{equation}
Based on Eqn. \eqref{recur-bary0}  and  the initialization value    $y_i(0)= \tilde{g}_i(x_i(0))$,   one can recursively show that
\begin{align}\label{recur-bary1}
 \bar{y}(k) = {1\over n} \sum_{i=1}^n \tilde{g}_i(x_i(k)) , \quad \forall k\geq 0.
\end{align}

{\em Step 1:  We first give an upper bound on $\| \bar{x}(k+1) -x^*\| $.}   From  Eqn. \eqref{recur-barx1}, it follows that
\begin{align}\label{recur-barx2}
&\| \bar{x}(k+1) -x^*\| = \| \bar{x}(k) -\alpha\bar{y}(k)-x^*\| \notag
\\& {(a)\atop = }\Big \| \bar{x}(k) -x^* - {\alpha \over n} \sum_{i=1}^n  \nabla f_i(\bar{x}(k)) + {\alpha \over n} \sum_{i=1}^n  \nabla f_i(\bar{x}(k))  \notag
\\&
\quad- {\alpha \over n} \sum_{i=1}^n  \nabla f_i(x_i(k))+  { \alpha\over n} \sum_{i=1}^n  \nabla f_i(x_i(k))  -\alpha \bar{y}(k)\Big \|
 \notag \\& {(b)\atop \leq  }  \left  \| \bar{x}(k) -\alpha  \nabla F(\bar{x}(k)) -x^* \right\|
 \notag\\&\quad +{ \alpha \over n} \left \|  \sum_{i=1}^n  \nabla f_i(\bar{x}(k)) -  \sum_{i=1}^n  \nabla f_i(x_i(k))\right\|
\notag\\& \quad +  \alpha \left\|  {  1\over n} \sum_{i=1}^n  \nabla f_i(x_i(k)) - \bar{y}(k)\right \|
\notag \\&  {(c)\atop \leq  }   \left  \| \bar{x}(k) -\alpha  \nabla F(\bar{x}(k)) -x^* \right\|  +    {\alpha \over n} \sum_{i=1}^n\left\| w_i(k)\right \|
\notag \\& \quad + { \alpha  L  \over n} \sum_{i=1}^n \| x_i(k) -  \bar{x}(k)  \| ,
\end{align}
where in (a) we added and subtracted $  {\alpha\over n} \sum_{i=1}^n  \nabla f_i(\bar{x}(k))$ and $  { \alpha\over n} \sum_{i=1}^n  \nabla f_i(x_i(k)) $,  in (b) we used the triangle
inequality and $F(x )= {1\over n } \sum_{i=1}^n f_i(x)$,  and  (c) is obtained by  using Eqn. \eqref{def-w},
Eqn. \eqref{recur-bary1}, and Assumption \ref{ass-fun}(ii).
We now  introduce an inequality from \cite[Eqn. (2.1.24)]{nesterov2013introductory} on   the
$\eta$-strongly convex and $L$-smooth   function  $f(x)$:
\begin{align}\label{inequ-f}
 & (x-y)^T(\nabla f(x) -\nabla f(y)) \geq {\eta L \over \eta+L} \|x-y\|^2   \notag \\&
+{1\over \eta+L}  \|\nabla f(x)- \nabla f(y)\|^2,\quad \forall x,y\in \mathbb{R}^d.
\end{align}
By  $\alpha  \in (0, {2\over \eta+L}]$, we have that  ${2\over \alpha}-\eta\geq L.$ Define $L' \triangleq {2\over \alpha}-\eta$.
 From Assumption \ref{ass-fun} it follows that    the    function $F(x)={1\over n} \sum_{i=1}^n f_i(x)$ is  $\eta$-strongly convex and $L'$-smooth.
Thus, by applying Eqn.  \eqref{inequ-f} with $x=x(k)$ and $y=x^*,$ from  $ \nabla F(x^*)=0$  and ${2\over \eta+L'}=\alpha $   it follows that
\begin{align*}
\begin{array}{ll}
 & \left \| \bar{x}(k) -\alpha  \nabla F(\bar{x}(k)) -x^* \right\|^2   \\&   =    \left \| \bar{x}(k)  -x^* \right\|^2+ \alpha^2 \left \|   \nabla F(\bar{x}(k))  \right\|^2
\\&  -2\alpha  (\bar{x}(k)   -x^*)^T (\nabla F(\bar{x}(k)) -\nabla F(x^*))
\\& \leq  \left \| \bar{x}(k)  -x^* \right\|^2 +\alpha^2 \left \|   \nabla F(\bar{x}(k))  \right\|^2
\\& - 2 \alpha \left({\eta L' \over \eta+L'} \|x_k-x^*\|^2+{1\over \eta+L'}  \|\nabla F(x_k) \|^2\right)
\\&  \leq  \Big(1-{2 \alpha\eta L' \over \eta+L'} \Big) \left \| \bar{x}(k)  -x^* \right\|^2
 \\& - \alpha \Big({2\over \eta+L'}-\alpha \Big)  \|\nabla F(x_k) \|^2
\\& \leq   \left(1-{2 \alpha\eta L' \over \eta+L'} \right) \left \| \bar{x}(k)  -x^* \right\|^2
\\& = \left(1- \alpha^2\eta L' \right)  \left \| \bar{x}(k)  -x^* \right\|^2=\left(1- \alpha\eta \right)^2  \left \| \bar{x}(k)  -x^* \right\|^2.
\end{array}\end{align*} Then we can bound the first term  of  Eqn. \eqref{recur-barx2}  by
 \begin{align}\label{bd-xstar}
& \left \| \bar{x}(k) -\alpha  \nabla F(\bar{x}(k)) -x^* \right\|  \leq \theta \left \| \bar{x}(k)  -x^* \right\|
\end{align}
  with $\theta \triangleq  1- \alpha\eta.$
Therefore, by  plugging Eqn. \eqref{bd-xstar}  into Eqn. \eqref{recur-barx2} and using the    following relation
\begin{align}\label{bd-me}
\sum_{i=1}^n \| e_i\|_2 \leq \sqrt{n} \| (e_1^T,\cdots, e_n^T)^T\|_2,
\end{align}
we make further modifications to Eqn. \eqref{recur-barx2}    as follows:
\begin{align}
&\| \bar{x}(k+1) -x^*\|   \leq   \theta \| \bar{x}(k)-x^*\|  \notag\\& + \alpha {L \over \sqrt{n}}  \| x(k) - (\mathbf{1} \otimes I_d)\bar{x}(k)  \| +    {\alpha \over n} \sum_{i=1}^n\left\| w_i(k)\right \| . \label{recur-barx21}
\end{align}

{\em Step 2: We give a bound on $ \| x(k+1) - (\mathbf{1} \otimes I_d)\bar{x}(k+1)  \|.$} Because  $A\mathbf{1}=\mathbf{1}$ and the   spectral radius     $\sigma_A$   of    $A-\mathbf{1}\mathbf{1}^T/n$  satisfies  $\sigma_A\in (0,1),$
for any $x\in \mathbb{R}^{nd}$ we have,
  \begin{equation}\label{bd-consensus}
\begin{split}
&  \| (A  \otimes I_d ) x-   (\mathbf{1} \otimes I_d)  \bar{x}\| \\&  = \| (A-\mathbf{1}\mathbf{1}^T/n)\otimes I_d  \big(x- (\mathbf{1} \otimes I_d)  \bar{x} \big)  \|  \\
&   \leq   \| A-\mathbf{1}\mathbf{1}^T/n \| \|   x- (\mathbf{1} \otimes I_d)  \bar{x} \|  \\&  \leq     \sigma_A \|   x - (\mathbf{1} \otimes I_d)  \bar{x}  \| ,
\end{split}
\end{equation}
where  $\bar{x}={1\over n} \sum_{i=1}^n x_i$.
This combined with  \eqref{push2},  \eqref{recur-barx1}, and  the triangle  inequality produces the following
\begin{equation}\label{recur-xerr}
\begin{array}{ll}
 & \| x(k+1) - (\mathbf{1} \otimes I_d)\bar{x}(k+1)  \|   \\& = \| (A  \otimes I_d ) x(k)-\alpha y(k)  - (\mathbf{1} \otimes I_d) (\bar{x}(k) -\alpha\bar{y}(k))\|     \\&  \leq   \|
  (A  \otimes I_d ) x(k)- (\mathbf{1} \otimes I_d) \bar{x}(k)\|  \\&+\| \alpha y(k)  -\alpha  (\mathbf{1} \otimes I_d)  \bar{y}(k) \|    \\&  \leq     \sigma_A \|   x(k)- (\mathbf{1} \otimes I_d)  \bar{x}(k)  \|
 + \alpha \|   y(k)- (\mathbf{1} \otimes I_d)  \bar{y}(k)\| .
\end{array}
\end{equation}

{\em Step 3: We give a bound on $\| y(k) - (\mathbf{1} \otimes I_d)  \bar{y}(k)\|$.}  From  Eqns. \eqref{def-w}, \eqref{def-fw}, and \eqref{recur-bary0} it   follows that
\begin{align*} \begin{array}{ll} &\bar{y}(k+1)-\bar{y}(k)={1\over n}\sum_{i=1}^n \big(w_i(k+1)+\nabla f_i(x_i(k+1)) \big)
\\&-{1\over n}\sum_{i=1}^n  \big(w_i(k)+\nabla f_i(x_i(k)) \big)
\\& = \left({1\over n}  \mathbf{1}^T\otimes I_d  \right) \big (\nabla (k+1)+w(k+1)-\nabla (k)-w(k ) \big).
\end{array}
\end{align*}
Then by using \eqref{pull2}, \eqref{bd-consensus},  $ \left \| I_n -{1\over n} \mathbf{1} \mathbf{1}^T \right \|  \leq 1$,   and  the triangle  inequality, we may obtain the following:
\begin{align} \label{recur-yerr}
 &  \| y(k+1) - (\mathbf{1} \otimes I_d)  \bar{y}(k+1)\|  \notag
\\&=  \big  \| (A  \otimes I_d )  y(k) +\nabla (k+1)+w(k+1)-\nabla (k)-w(k ) \notag
\\ &  - (\mathbf{1} \otimes I_d) \bar{y}(k) +  (\mathbf{1} \otimes I_d)  (\bar{y}(k)-\bar{y}(k+1))   \big\| \notag
\\&\leq \| (A  \otimes I_d )  y(k)  - (\mathbf{1} \otimes I_d) \bar{y}(k)  \|  +\notag
 \\&  \big\| \big(I_n -{1\over n} \mathbf{1} \mathbf{1}^T \big) \otimes I_d (\nabla (k+1)+w(k+1)-\nabla (k)-w(k )  ) \big\|  \notag
\\&  \leq  \sigma_A \|   y(k)- (\mathbf{1} \otimes I_d)  \bar{y}(k)  \|   + \notag
\\&  \big \| I_n -{1\over n} \mathbf{1} \mathbf{1}^T \big \|  \big ( \|  \nabla (k+1)-\nabla (k) \|+ \| w(k+1)-w(k ) \|  \big )  \notag
\\&  \leq  \sigma_A \|   y(k)- (\mathbf{1} \otimes I_d)  \bar{y}(k)  \|\notag
   \\& +   L \|x(k+1)-x(k)\| +  \| w(k+1)-w(k ) \| ,
\end{align}
where in the last inequality we used   the Lipschitz continuity of  $\nabla f_i$  (Assumption \ref{ass-fun}(ii)) and the definition of $\nabla (k) $ in \eqref{def-w}.
\begin{align*}
&\|  \nabla (k+1)-\nabla (k) \|
\\& =\sqrt{\sum_{i=1}^n \|\nabla f_i(x_i(k+1))-\nabla f_i(x_i(k)) \|^2}
\\& \leq \sqrt{\sum_{i=1}^nL^2  \| x_i(k+1)- x_i(k) \|^2}= L  \| x(k+1)- x(k) \| .
\end{align*}
We  then give an estimate for the upper bound of $ \|x(k+1)-x(k)\| .$ From \eqref{push2} and $A\mathbf{1} =\mathbf{1} $ it follows that
\begin{align}
&  \| x (k+1)-x(k) \|   = \| (A  \otimes I_d )x(k)-\alpha y(k)   -x(k)\| \notag \\& =  \| (A-I_n)  \otimes I_d \big(x(k)-(\mathbf{1} \otimes I_d) \bar{x}(k)\big) \notag
\\& -\alpha \big( y(k)- (\mathbf{1} \otimes I_d)  \bar{y}(k) \big)  -\alpha  (\mathbf{1} \otimes I_d)  \bar{y}(k)  \|  \notag
\\& \leq   \|  A-I_n \|  \big \| x(k)-(\mathbf{1} \otimes I_d) \bar{x}(k)\big\|  \notag
\\&+\alpha  \|   y(k)- (\mathbf{1} \otimes I_d)  \bar{y}(k)  \|   +\alpha \sqrt{ n} \|    \bar{y}(k)  \| , \label{recur-xdiff}
 \end{align}
where in the last inequality  we used the triangle inequality   and $\|(\mathbf{1} \otimes I_d)  \bar{y}(k)  \|  =\sqrt{n} \|    \bar{y}(k)  \| . $
Then  by substituting Eqn. \eqref{recur-xdiff} into  Eqn. \eqref{recur-yerr}  there holds
\begin{equation}  \label{recur-yerr2}
\begin{array}{ll}
  &\| y(k+1) - (\mathbf{1} \otimes I_d)  \bar{y}(k+1)\|
\\&  \leq ( \sigma_A+\alpha L)  \|   y(k)- (\mathbf{1} \otimes I_d)  \bar{y}(k)  \|
\\&+   L    \|  A-I_n \|  \big \| x(k)-(\mathbf{1} \otimes I_d) \bar{x}(k)\big\|
\\&     +\alpha \sqrt{ n} L \|    \bar{y}(k)  \| +  \| w(k+1)-w(k ) \| .
\end{array}
\end{equation}
Next, we   provide an  upper bound on $\|    \bar{y}(k)  \| .$ By    using $\sum_{i=1}^n \nabla f_i(x^*)=0, $  Eqns. \eqref{def-w} and \eqref{recur-bary1}, we obtain that
\begin{align*}
&\|    \bar{y}(k)  \|=\big\|  {1\over n} \sum_{i=1}^n (  w_i(k) + \nabla f_i(x_i(k)) )\big \|  \notag
\\&  \leq \big \|  {1\over n} \sum_{i=1}^n w_i(k) \big \|  + \big \|  {1\over n} \sum_{i=1}^n  (\nabla f_i(x_i(k))-\nabla f_i(x^*))  \big \|   \notag
\\&  {(a) \atop \leq} {1\over n} \sum_{i=1}^n  \left \|  w_i(k) \right \|  + { L\over n} \sum_{i=1}^n  \| x_i(k) - x^*\|
\\&  {(b) \atop \leq} {1\over n} \sum_{i=1}^n \left \|   w_i(k) \right \|  + { L\over \sqrt{n}}   \| x(k) - (\mathbf{1} \otimes I_d)   x^*\|
\notag \\ &
{(c) \atop  \leq}{1\over n} \sum_{i=1}^n \left \|   w_i(k) \right \|    + { L\over \sqrt{n}}   \| x(k) - (\mathbf{1} \otimes I_d)   \bar{x}(k)\|
\\&\quad+{ L\over \sqrt{n}}   \| (\mathbf{1} \otimes I_d)   \bar{x}(k) - (\mathbf{1} \otimes I_d)   x^*\|  \\ &
  \leq {1\over n} \sum_{i=1}^n \left \|   w_i(k) \right \|  + { L\over \sqrt{n}}   \| x(k) - (\mathbf{1} \otimes I_d)   \bar{x}(k)\| \\&
  \quad +  L   \|   \bar{x}(k) -   x^*\|,
\end{align*}
where in (a)   we used  Assumption \ref{ass-fun}(ii),   in (b) we  utilized Eqn.  \eqref{bd-me},  and in (c)  we added and  subtracted  the term $(\mathbf{1} \otimes I_d)   \bar{x}(k) $ and applied the triangle inequality.
%, and $ \| (\mathbf{1} \otimes I_d)   \bar{x}(k) - (\mathbf{1} \otimes I_d)   x^*\|=\sqrt{n} \|   \bar{x}(k) -   x^*\|$ is used in the last inequality.
  This combined with  \eqref{recur-yerr2}  produces
\begin{align}  \label{recur-yerr3}
  &  \| y(k+1) - (\mathbf{1} \otimes I_d)  \bar{y}(k+1)\|   \notag
\\&  \leq ( \sigma_A+\alpha L)  \|   y(k)- (\mathbf{1} \otimes I_d)  \bar{y}(k)  \|    + \alpha \sqrt{n} L^2  \|   \bar{x}(k) -   x^*\|  \notag
  \\& +  \left(   L    \|  A-I_n \|  + \alpha L^2\right)\big \| x(k)-(\mathbf{1} \otimes I_d) \bar{x}(k)\big\|   \notag
\\& +  \| w(k+1)-w(k ) \|+{\alpha   L \over \sqrt{n}} \sum_{i=1}^n \left \|   w_i(k) \right \| .
\end{align}
{\em  Step 4: Obtain a system of  inequalities.} By the definition of $z(k)$ as in  \eqref{def-z}, and by
 combining  Eqns. \eqref{recur-barx21}, \eqref{recur-xerr},   and \eqref{recur-yerr3},  we obtain that
\begin{align*}
z(k+1) &\leq  \begin{pmatrix}
 \theta &  \alpha {L \over \sqrt{n}}  & 0\\ 0&   \sigma_A & \alpha  \\
 \alpha \sqrt{ n} L^2 &L    \|  A-I_n \|  + \alpha L^2  &\sigma_A+\alpha L \end{pmatrix} z(k)
  \\& +\begin{pmatrix}
 {\alpha \over n} \sum_{i=1}^n\left\| w_i(k)\right \|\\
0 \\ \| w(k+1)-w(k ) \|  +{\alpha   L \over \sqrt{n}} \sum_{i=1}^n \left \|   w_i(k) \right \|
\end{pmatrix} .
\end{align*}
Then by recalling  $\theta=1-\alpha L$,  we prove the lemma.
\hfill $\Box$
\bibliographystyle{IEEEtran}
\bibliography{DSO}

 \end{document}